\documentclass[12pt,draftcls,onecolumn]{IEEEtranTOAPPEAR}

\usepackage{amssymb,latexsym,color,amsmath,pifont,dsfont,epsfig,mathtools}
\usepackage{graphicx}
\usepackage{subfig}
\usepackage{setspace,cite}
\usepackage{algorithm}
\usepackage{algorithmic}

\usepackage{hyperref}
\hypersetup{colorlinks=false,pdfborder={0 0 0},bookmarks=false}
\newcommand{\non}{\nonumber}

\newcommand{\ba}{\begin{array}}
\newcommand{\ea}{\end{array}}
\newcommand{\be}{\begin{equation}}
\newcommand{\ee}{\end{equation}}
\newcommand{\bes}{\begin{equation*}}
\newcommand{\ees}{\end{equation*}}
\newcommand{\bea}{\begin{eqnarray}}
\newcommand{\eea}{\end{eqnarray}}
\newcommand{\beas}{\begin{eqnarray*}}
\newcommand{\eeas}{\end{eqnarray*}}
\newcommand{\ben}{\begin{enumerate}}
\newcommand{\een}{\end{enumerate}}
\newcommand{\bc}{\begin{center}}
\newcommand{\ec}{\end{center}}	
\newcommand{\bi}{\begin{itemize}}
\newcommand{\ei}{\end{itemize}}
\newcommand{\bex}{\begin{example}}
\newcommand{\eex}{\end{example}}
\newcommand{\br}{\begin{remark}}
\newcommand{\er}{\end{remark}}
\newcommand{\bd}{\begin{definition}}
\newcommand{\ed}{\end{definition}}
\newcommand{\bt}{\begin{theorem}}
\newcommand{\et}{\end{theorem}}
\newcommand{\bl}{\begin{lemma}}
\newcommand{\el}{\end{lemma}}
\newcommand{\bp}{\begin{proof}}
\newcommand{\ep}{\end{proof}}

\newcommand{\beqns}{\begin{eqnarray*}}	%old; use \beas
\newcommand{\eeqns}{\end{eqnarray*}}		%old; use \eeas
\newcommand{\bct}{\begin{center}}		%old; use \bc
\newcommand{\ect}{\end{center}}			%old; use \ec
\newcommand{\benum}{\begin{enumerate}}	%old; use \ben
\newcommand{\eenum}{\end{enumerate}}		%old; use \een

\newcommand{\veps}{\varepsilon}

\newcommand{\enma}[1]   {\ensuremath{#1}}

\newcommand{\trace}     {\enma{\mathrm{trace}}}

\newcommand{\eS}{{\cal S}}

\newtheorem{theorem}{Theorem}
\newtheorem{lemma}[theorem]{Lemma}

\newtheorem{proposition}[theorem]{Proposition}

\newtheorem{definition}{Definition}

\newtheorem{remark}{Remark}
\newtheorem{example}{Example}

\newcommand{\matbegin}{
        \left[
}
\newcommand{\matend}{
        \right]
}

\newcommand{\tbo}[2]{
  \matbegin \begin{array}{c}
       #1 \\ #2
       \end{array} \matend }

\newcommand{\cL}{{\cal L}}

\newcommand{\gs}{G_{ij}^\star}
\newcommand{\vij}{V_{ij}}
\newcommand{\gij}{G_{ij}}

\newcommand{\htwo}{{\cal H}_2}
\newcommand{\card}{\enma{{\bf card}}}
\newcommand{\mre}{\mathrm{e}}
\newcommand{\mrd}{\mathrm{d}}

\IEEEoverridecommandlockouts
\title{\LARGE \bf
Design of Optimal Sparse Feedback Gains via \\[-0.3cm]
the Alternating Direction Method of Multipliers}

\author{Fu Lin, Makan Fardad, and Mihailo R. Jovanovi\'c
\thanks{Financial support from the National Science Foundation under
awards CMMI-0927509 and CMMI-0927720 and under CAREER Award CMMI-0644793
is gratefully acknowledged.}
\thanks{F.\ Lin and M.\ R.\ Jovanovi\'c are with the Department of
Electrical and Computer Engineering, University of Minnesota,
Minneapolis, MN 55455. M.\ Fardad is with the Department of Electrical Engineering
and Computer Science, Syracuse University, NY 13244. E-mails: fu@umn.edu, makan@syr.edu, mihailo@umn.edu.}}

\begin{document}
\maketitle

\vspace*{-1.75cm}

    \begin{abstract}
    \vspace*{-2ex}
We design sparse and block sparse feedback gains that minimize the {variance amplification (i.e., the} $\htwo$ norm) of distributed systems. Our approach consists of two steps. First, we identify sparsity patterns of feedback gains by incorporating sparsity-promoting penalty functions into the {optimal control} problem, where the added terms penalize the number of communication links in the distributed controller. Second, we optimize feedback gains subject to structural constraints determined by the identified sparsity patterns. In the first step, the sparsity structure of feedback gains is identified using the alternating direction method of multipliers, which is a powerful algorithm well-suited to large optimization problems. This method alternates between {promoting} the sparsity {of the controller} and optimizing the closed-loop {performance}, which allows us to exploit the structure of the corresponding objective functions. In particular, we take advantage of the separability of the sparsity-promoting penalty functions to decompose the minimization problem into sub-problems that can be solved analytically. Several examples are provided to illustrate the effectiveness of the developed approach.
    \end{abstract}

    \vspace*{-2ex}
\begin{keywords}
    \vspace*{-2ex}
Alternating direction method of multipliers, communication architectures, continuation methods, $\ell_1$ minimization, optimization, separable penalty functions, sparsity-promoting optimal control, structured distributed design.
\end{keywords}

    \vspace*{-2ex}
\section{Introduction}
    \label{sec.intro}

	We develop methods for the design of {\em sparse\/} and {\em block sparse\/} feedback gains that minimize the variance amplification of distributed systems. Our approach consists of two steps. The first step, which can be viewed as a {\em structure identification step\/}, is aimed at finding sparsity patterns that strike a balance between the $\htwo$ performance and the sparsity of the controller. This is achieved by incorporating {\em sparsity-promoting\/} penalty functions into the optimal control problem, where the added sparsity-promoting terms penalize the number of communication links. We consider several sparsity-promoting penalty functions including the cardinality function and its convex relaxations. In the absence of sparsity-promoting terms, the solution to the standard $\htwo$ problem results in centralized controllers with dense feedback gains. By gradually increasing the weight on the sparsity-promoting penalty terms, the optimal feedback gain moves along a parameterized solution path from the centralized to the sparse gain of interest. This weight is increased until the desired balance between performance and sparsity is achieved. In the second step, in order to improve the $\htwo$ performance of the structured controller, we solve an optimal control problem subject to the feedback gain belonging to the identified structure.

	We demonstrate that the {\em alternating direction method of multipliers\/}~(ADMM)~\cite{boyparchupeleck11} provides an effective tool for the design of sparse distributed controllers whose performance is comparable to the performance of the optimal centralized controller. This method alternates between {promoting} the sparsity of the feedback gain matrix and optimizing the closed-loop $\htwo$ norm. The advantage of this alternating mechanism is threefold. First, it provides a flexible framework for incorporation of different penalty functions that promote sparsity or block sparsity. Second, it allows us to exploit the {\em separability\/} of the sparsity-promoting penalty functions and to decompose the corresponding optimization problems into sub-problems that can be solved {\em analytically\/}. {These analytical results are immediately applicable to other distributed control problems where sparsity is desired.} Finally, it facilitates the use of descent algorithms for $\htwo$ optimization, in which a descent direction can be formed by solving two Lyapunov equations and one Sylvester equation.

	The $\ell_1$ norm is widely used as a proxy for cardinality minimization in applied statistics, in sparse signal processing, and in machine learning; see~\cite{cantao06,hastibfri09,boyparchupeleck11,canwakboy08}. In the controls community, recent work inspired by similar ideas includes~\cite{schebeall11,schlilamall11,zelschall13}. In~\cite{schebeall11}, an $\ell_0$ induced gain was introduced to quantify the sparsity of the impulse response of a discrete-time system. In~\cite{schlilamall11}, the weighted $\ell_1$ framework was used to design structured dynamic output feedback controllers subject to a given ${\cal H}_\infty$ performance. In~\cite{zelschall13}, an $\ell_1$ relaxation method was employed for the problem of adding a fixed number of edges to a consensus network.

	Our presentation is organized as follows. We formulate the sparsity-promoting optimal control problem and compare several sparsity-promoting penalty functions in Section~\ref{sec.SP}. We present the ADMM algorithm, emphasize the separability of the penalty functions, and provide the analytical solutions to the sub-problems for both sparse and block sparse minimization problems in Section~\ref{sec.ADMM}. Several examples are provided in Section~\ref{sec.example} to demonstrate the effectiveness of the developed approach. We conclude with a summary of our contributions in Section~\ref{sec.conclusion}.

    \vspace*{-2ex}
\section{Sparsity-promoting Optimal Control Problem}
    \label{sec.SP}

Consider the following control problem
    \begin{align}
    \non
    \dot{x}
    \;&=\;
    A \, x \,+\, B_1 \, d \,+\, B_2 \, u
    \\[-0.25cm]
    z
    \;&=\;
    C \, x \,+\, D \, u
    \label{eq.ss}
    \\[-0.25cm]
    u
    \;&=\;
    - \,
    F \, x
    \non
    \end{align}
where $d$ and $u$ are the disturbance and control inputs, $z$ is the performance output, $C = \big[\, Q^{1/2} ~~ 0 \,\big]^T$, and $D = \big[\, 0 ~~ R^{1/2} \,\big]^T$, {with standard assumptions that $(A,B_2)$ is stabilizable and $(A,Q^{1/2})$ is detectable.} The matrix $F$ is a state feedback gain,  $Q = Q^T \geq 0$ and $R = R^T > 0$ are the state and control performance weights, and the closed-loop system is given by
    \be
    \ba{rcl}
    \dot{x}
    & = &
    (A \,-\, B_2 F) \, x \,+\, B_1 \, d
    \\[0.1cm]
    z
    & = &
    \left[
    \ba{c}
    \!\! Q^{1/2} \!\! \\
    \!\! -\,R^{1/2} F \!\!
    \ea
    \right]
    x.
    \ea
    \label{CL}
    \ee

The design of the optimal state feedback gain $F$, subject to structural constraints that dictate its zero entries, was recently considered by the authors in~\cite{farlinjovCDC09,linfarjovTAC11al}. Let the subspace $\eS$ embody these constraints and let us assume that there exists a stabilizing $F \in \eS$. References~\cite{farlinjovCDC09,linfarjovTAC11al} then search for $F \in \eS$ that minimizes the $\htwo$ norm of the transfer function from $d$ to $z$,
    \be
    \ba{ll}
    \!\!\!\!
    \text{minimize}
    & J(F)
    \\
    \!\!\!\!
    \text{subject to}
    &
    F \in \eS
    \ea
    \tag*{(SH2)}
    \label{SH2}
    \ee
where
    \be
    \label{eq.J}
    J(F)
    \, = \,
    \left\{
    \ba{rl}
    \trace \,
    \big(
    B_1^T P(F) B_1
    \big),
    &
    \mbox{$F$ stabilizing}
    \\
    \infty,
    &
    \mbox{otherwise}.
    \ea
    \right.
    \ee
The matrix $P (F)$ in~(\ref{eq.J}) denotes the closed-loop observability Gramian
    \be
        P(F)
        \; = \;
        \int_{0}^\infty
        \mre^{
        (A - B_2 F )^T t
        }
        \,
        (Q + F^T R F )
        \,
        \mre^{
        (A - B_2 F) t
        }
        \, \mrd t
    \ee
which can be obtained by solving the Lyapunov equation
    \be
    \label{eq.lyap}
    \left( A - B_2 F \right)^T P
    \,+\,
    P
    \left( A - B_2 F \right)
    \,=\,
    -\left( Q + F^T R F \right).
    \ee

While the communication architecture of the controller in~\ref{SH2} is {\em a priori\/} specified, in this note our emphasis shifts to identifying favorable communication structures without any prior assumptions on the sparsity patterns of the matrix $F$. We propose an optimization framework in which the sparsity of the feedback gain is directly incorporated into the objective function.

Consider the following optimization problem
    \be
    \text{minimize}
    \;\;
    J(F)
    \,+\,
    \gamma
    \,
    g_0(F)
    \label{card}
    \ee
where
    \be
    \label{eq.card}
    g_0(F)
    \, = \,
    \card
    \left( F \right)
    \ee
denotes the cardinality function, i.e., the {\em number of nonzero elements\/} of a matrix. In contrast to problem~\ref{SH2}, no structural constraint is imposed on $F$; instead, our goal is to promote sparsity of the feedback gain by incorporating the cardinality function into the optimization problem. The positive scalar $\gamma$ characterizes our emphasis on the sparsity of $F$; a larger $\gamma$ encourages a sparser $F$, while $\gamma = 0$ renders a centralized gain that is the solution of the standard LQR problem. For $\gamma = 0$, the solution to~(\ref{card}) is given by $F_c = R^{-1} B_2^T P$, where $P$ is the unique positive definite solution of the algebraic Riccati equation,
	$
    A^T P
    +
    P A
    +
    Q
    -
    P B_2 R^{-1} B_2^T P
    =
    0.
    	$

    \vspace*{-2ex}
\subsection{Sparsity-promoting penalty functions}
    \label{sec.penalty}

Problem~(\ref{card}) is a combinatorial optimization problem whose solution usually requires an intractable combinatorial search. In optimization problems where sparsity is desired, the cardinality function is typically replaced by the $\ell_1$ norm of the optimization variable~\cite[Chapter 6]{boyvan04},
    \be
    \label{eq.l1}
    g_1(F)
    \, = \,
    \| F \|_{\ell_1}
    \, = \,
    \sum_{i, \, j}
    | F_{ij} |.
    \ee
Recently, a {\em weighted\/} $\ell_1$ norm was used to enhance sparsity in signal recovery~\cite{canwakboy08},
    \be
    \label{eq.wl1}
    g_2(F)
    \, = \,
    \sum_{i, \, j}
    W_{ij}
    |F_{ij}|
    \ee
where $W_{ij}$ are non-negative weights. If $W_{ij}$'s are chosen to be inversely proportional to the magnitude of $F_{ij}$, i.e.,
    $
    \{ W_{ij} = 1/|F_{ij}|$,
    $F_{ij} \neq 0$;
    $W_{ij} = 1/\veps$,
    $F_{ij} = 0$,
    $0 < \veps \ll 1 \}$,
then the weighted $\ell_1$ norm and the cardinality function of $F$ coincide,
    $
    \sum_{i, \, j} W_{ij} \, |F_{ij}|
    =
    \card
    \left( F \right).
    $
This scheme for the weights, however, cannot be implemented, since the weights depend on the unknown feedback gain. A reweighted algorithm that solves a sequence of weighted $\ell_1$ optimization problems in which the weights are determined by the solution of the weighted $\ell_1$ problem in the previous iteration was proposed in~\cite{lobfazboy07,canwakboy08}. This reweighted scheme was recently employed by the authors to design sparse feedback gains \mbox{for a class of distributed systems~\cite{farlinjovACC11,farlinjovACC12}.}

Both the $\ell_1$ norm and its weighted version are {\em convex\/} relaxations of the cardinality function. On the other hand, we also examine utility of the {\em nonconvex\/} {\em sum-of-logs\/} function as a more aggressive means for promoting sparsity~\cite{canwakboy08}
    \be
    \label{eq.logsum}
    g_3(F)
    \, = \,
    \sum_{i, \, j}
    \log
    \left(
    1
    \, + \,
    |F_{ij}| / \veps
    \right),
    ~~~
    0 \, < \, \veps \, \ll \, 1.
    \ee

    \remark
    \label{rem.g_group}
{Design of feedback gains that have block sparse structure can be achieved by promoting sparsity at the level of the {\em submatrices\/} instead of at the level of the individual elements. Let the feedback gain $F$ be partitioned into submatrices $F_{ij} \in \mathbb{R}^{m_i \times n_j}$ that need not have the same size. The weighted $\ell_1$ norm and the sum-of-logs can be generalized to matrix blocks by replacing the absolute value of $F_{ij}$ in~(\ref{eq.wl1}) and~(\ref{eq.logsum}) by the Frobenius norm $\| \cdot \|_F$ of $F_{ij}$. Similarly, the cardinality function~(\ref{eq.card}) should be replaced by
    $
    \sum_{i,j}
    \card
    \left(
    \| F_{ij} \|_F
    \right),
    $
where $\| F_{ij} \|_F$ does not promote sparsity within the $F_{ij}$ block; \mbox{it instead promotes sparsity at the level of submatrices.}}

    \vspace*{-2ex}
\subsection{Sparsity-promoting optimal control problem}

Our approach to sparsity-promoting feedback design makes use of the above discussed penalty functions. In order to obtain state feedback gains that {strike a balance between the quadratic performance and the sparsity of the controller}, we consider the following optimal control problem
    \be
    \text{minimize}
    \;\;
    J(F)
    \, + \,
    \gamma
    \,
    g(F)
    \tag*{(SP)}
    \label{SP}
    \ee
	where $J$ is the square of the closed-loop $\htwo$ norm~(\ref{eq.J}) and $g$ is a sparsity-promoting penalty function, e.g., given by~(\ref{eq.card}), (\ref{eq.l1}), (\ref{eq.wl1}), or (\ref{eq.logsum}). When the cardinality function in~(\ref{eq.card}) is replaced by (\ref{eq.l1}), (\ref{eq.wl1}), or (\ref{eq.logsum}), problem~\ref{SP} can be viewed as a relaxation of the combinatorial problem~(\ref{card})-(\ref{eq.card}), obtained by approximating the cardinality function \mbox{with the corresponding penalty functions $g$.}

	As the parameter $\gamma$ varies over $[0, +\infty)$, the solution of~\ref{SP} traces the trade-off path between the $\htwo$ performance $J$ and the feedback gain sparsity $g$. When $\gamma = 0$, the solution is the centralized feedback gain $F_c$. We then slightly increase $\gamma$ and employ an iterative algorithm -- the alternating direction method of multipliers (ADMM) -- initialized by the optimal feedback matrix at the previous $\gamma$. The solution of~\ref{SP} becomes sparser as $\gamma$ increases. After a desired level of sparsity is achieved, we fix the sparsity structure and find the optimal structured feedback gain by solving the structured $\htwo$ problem~\ref{SH2}.

Since the set of stabilizing feedback gains is in general not convex~\cite{perger94} and since the matrix exponential is not necessarily a convex function of its argument~\cite{boyvan04}, $J$ need not be a convex function of $F$. This makes it difficult to establish convergence to the global minimum of~\ref{SP}. Even in problems for which we cannot establish the convexity of $J (F)$, our extensive computational experiments suggest that the algorithms developed in Section~\ref{sec.ADMM} provide an effective means for attaining a desired trade-off between the \mbox{$\htwo$ performance and the sparsity of the controller.}

    \vspace*{-2ex}
\section{Identification of Sparsity-patterns via ADMM}
    \label{sec.ADMM}

Consider the following constrained optimization problem
    \be
    \label{eq.constr}
    \ba{ll}
    \!\!\!\!
    \text{minimize}
    & J(F)
    \, + \,
    \gamma
    \,
    g(G)
    \\
    \!\!\!\!
    \text{subject to}
    &
    F
    \,-\,
    G
    \,=\,
    0
    \ea
    \ee
which is clearly equivalent to the problem~\ref{SP}. The augmented Lagrangian associated with the constrained problem~(\ref{eq.constr}) is given by
    \[
    \ba{rcl}
    \cL_\rho(F,G,\Lambda)
    &=&
    J(F)
    \,+\,
    \gamma
    \,
    g(G)
    \,+\,
    \trace
    \left(
    \Lambda^T ( F \,-\, G )
    \right)
    \, + \,
    (\rho / 2)
    \| F \,-\, G \|_F^2
    \ea
    \]
where $\Lambda$ is the dual variable (i.e., the {\em Lagrange multiplier\/}), $\rho$ is a positive scalar, and $\| \cdot \|_{F}$ is the Frobenius norm. By introducing an additional variable $G$ and an additional constraint $F - G = 0$, we have simplified the problem~\ref{SP} by decoupling the objective function into two parts that depend on two different variables. As discussed below, this allows us to exploit the structures of $J$ and $g$.

In order to find a minimizer of the constrained problem~(\ref{eq.constr}), the ADMM algorithm uses a sequence of iterations
    \begin{subequations}
    \label{eq.ADMM}
    \begin{align}
    \label{eq.F_update}
    F^{k+1}
    \,&:=\,
    \underset{F}{\operatorname{arg \, min}}
    \;
    \cL_\rho
    (F,G^k,\Lambda^k)
    \\[-0.2cm]
    \label{eq.G_update}
    G^{k+1}
    \,&:=\,
    \underset{G}{\operatorname{arg \, min}}
    \;
    \cL_\rho
    (F^{k+1},G,\Lambda^k)
    \\[-0.2cm]
    \label{eq.lambda_update}
    \Lambda^{k+1}
    \,&:=\,
    \Lambda^{k}
    \,+\,
    \rho ( F^{k+1} \,-\, G^{k+1})
    \end{align}
    \end{subequations}
until
    $
    \| F^{k+1} - G^{k+1} \|_F
    \leq
    \epsilon
    $
and
    $
    \| G^{k+1} - G^{k} \|_F
    \leq
    \epsilon.
    $
In contrast to the {\em method of multipliers\/}~\cite{boyparchupeleck11}, in which $F$ and $G$ are {\em minimized jointly\/},
   $
    (
    F^{k+1},
    G^{k+1}
    )
    :=
    \underset{F,\, G}{\operatorname{arg \, min}}
    \;
    \cL_\rho
    (F,G,\Lambda^k),
    $
ADMM consists of an $F$-minimization step~(\ref{eq.F_update}), a $G$-minimization step~(\ref{eq.G_update}), and a dual variable update step~(\ref{eq.lambda_update}). Note that the dual variable update~(\ref{eq.lambda_update}) uses a step-size equal to $\rho$, which guarantees {that one of the dual feasibility conditions is satisfied in each ADMM iteration; see~\cite[Section 3.3]{boyparchupeleck11}.}

ADMM brings two major benefits to the sparsity-promoting optimal control problem~\ref{SP}:
    \bi
    \item {\em Separability of $g$.\/} The penalty function $g$ is {\em separable\/} with respect to the {\em individual\/} elements of the matrix. In contrast, the closed-loop $\htwo$ norm cannot be decomposed into componentwise functions of the feedback gain. By separating $g$ and $J$ in the minimization of the augmented Lagrangian $\cL_\rho$, we can determine {\em analytically\/} the solution to the $G$-minimization problem via decomposition of~(\ref{eq.G_update}) into sub-problems \mbox{that only involve {\em scalar\/} variables.}

    \item {\em Differentiability of $J$.\/} The square of the closed-loop $\htwo$ norm $J$ is a {\em differentiable\/} function of $F$~\cite{linfarjovTAC11al}; this is in contrast to $g$ which is a {\em non-differentiable\/} function. By separating $g$ and $J$ in the minimization of the augmented Lagrangian $\cL_\rho$, we can utilize descent algorithms that rely on the differentiability of $J$ to solve the $F$-minimization problem~(\ref{eq.F_update}).
    \ei

{We next {provide} the analytical expressions for the solutions of the $G$-minimization problem~(\ref{eq.G_update}) in Section~\ref{sec.G_update}, describe a descent method to solve the $F$-minimization problem~(\ref{eq.F_update}) in Section~\ref{sec.F_update}, present Newton's method to solve the structured problem~\ref{SH2} in Section~\ref{sec.polishing}, and discuss the convergence of ADMM in Section~\ref{sec.convergence}.}

    \vspace*{-2ex}
\subsection{Separable solution to the $G$-minimization problem~{\em (\ref{eq.G_update})\/}}
    \label{sec.G_update}

The completion of squares with respect to $G$ in the augmented Lagrangian $\cL_\rho$ can be used to show that~(\ref{eq.G_update}) is equivalent to
    \be
    \label{eq.min_Gij}
    \text{minimize}
    ~~
    \phi(G)
    \, = \,
    \gamma
    \,
    g(G)
    \,+\,
    (\rho/2)
    \| G \,-\, V^{k} \|_F^2
    \ee
where
    $
    V^{k}
    =
    (1/\rho)
    \Lambda^k
    +
    F^{k+1}.
    $
To simplify notation, we drop the superscript $k$ in $V^{k}$ throughout this section. Since both $g$ and the square of the Frobenius norm can be written as a summation of componentwise functions of a matrix, we can decompose~(\ref{eq.min_Gij}) into sub-problems expressed in terms of the {\em individual\/} elements of $G$. For example, if $g$ is the weighted $\ell_1$ norm, then
    $
    \phi(G)
    =
    \sum_{i,j}
    \left(
    \gamma \, W_{ij} \, |G_{ij}| + (\rho/2) (G_{ij} - V_{ij})^2
    \right).
    $
This facilitates the conversion of (\ref{eq.min_Gij}) to minimization problems that only involve {\em scalar\/} variables $G_{ij}$. By doing so, the solution of~(\ref{eq.min_Gij}) can be determined {\em analytically\/} for the weighted $\ell_1$ norm, the sum-of-logs, and the cardinality function.

\subsubsection{Weighted $\ell_1$ norm}

The unique solution to~(\ref{eq.min_Gij}) is given by the {\em soft thresholding\/} operator~(e.g., see~\cite[Section 4.4.3]{boyparchupeleck11})
    \be
    \label{eq.shrinkage_wl1}
    G_{ij}^\star
    \, = \,
    \left\{
    \ba{ll}
    \left(
    1
    \, - \,
    {a}/{| \vij |}
    \right)
    \vij,
    &
    | \vij |
    \, > \,
    a
    \\
    0,
    &
    | \vij |
    \, \leq \,
    a
    \ea
    \right.
    \ee where $a = (\gamma/\rho) W_{ij}$. For given $V_{ij}$, $G_{ij}^\star$ is obtained by moving $V_{ij}$ towards zero with the amount $(\gamma/\rho) W_{ij}$. In particular, $G_{ij}^\star$ is set to zero if $|V_{ij}| \leq (\gamma/\rho) W_{ij}$, implying that a more aggressive scheme for driving $\gs$ to zero can be obtained by increasing $\gamma$ and $W_{ij}$ and by decreasing $\rho$.

\subsubsection{Cardinality function}
    \label{sec.card}

The unique solution to~(\ref{eq.min_Gij}) is given by the {\em truncation\/} operator
    \be
    \label{eq.shrinkage_card}
    G_{ij}^\star
    \, = \,
    \left\{
    \ba{ll}
    V_{ij},
    &
    |V_{ij}|
    \, > \,
    b
    \\
    0,
    &
    |V_{ij}|
    \, \leq \,
    b
    \ea
    \right.
    \ee
where $b = \sqrt{2 \gamma / \rho}$. For given $\vij$, $\gs$ is set to $\vij$ if $|\vij| > \sqrt{2 \gamma / \rho}$ and to zero if $|\vij| \leq \sqrt{2 \gamma / \rho}$.

\subsubsection{Sum-of-logs function}

As shown in~\cite{fu-phd12}, the solution to~(\ref{eq.min_Gij}) is given by
        \be
        \label{eq.shrinkage_logsum}
        G_{ij}^\star
        \,=\,
        \left\{
        \ba{ll}
        0,
        &
        \Delta
        \, \leq \,
        0
        ~{\rm or}~
        \{
        \Delta \, > \, 0
        ~{\rm and}~
        r^+
        \}
        \\
        r^+ \, V_{ij},
        &
        \Delta \, > \, 0
        ~{\rm and}~
        r^- \, \leq \, 0
        ~{\rm and}~
        0 \, < \,
        r^+
        \, \leq \, 1
        \\
        G^0,
        &
        \Delta \, > \, 0
        ~{\rm and}~
        0 \, \leq \,
        r^\pm
        \, \leq \, 1
        \ea
        \right.
        \ee
where
        \be
        \label{eq.rpm}
        \ba{rcl}
        \Delta
        & = &
        (|\vij| + \veps)^2 - 4(\gamma/\rho)
        \\
        r^{\pm}
        & = &
        \left(
        | V_{ij} |
        \, - \,
        \veps
        \, \pm \,
        \sqrt{\Delta}
        \right) / \left( {2 \, | V_{ij} |} \right)
        \ea
        \ee
and
        $
        G^0
        :=
        \operatorname{arg \, min}
        \,
        \{
        \phi_{ij}(r^+ \vij),
        \phi_{ij}(0)
        \}.
        $
For fixed $\rho$ and $\veps$, (\ref{eq.shrinkage_logsum}) is determined by the value of $\gamma$. For small $\gamma$, (\ref{eq.shrinkage_logsum}) resembles the soft thresholding operator and for large $\gamma$, it resembles the truncation operator.

    \remark
{In block sparse design, $g$ is determined by
    $
    \left\{
    \sum_{i,j}
    W_{ij}
    \| G_{ij} \|_F
    \right.
    $;
    $
    \sum_{i,j}
    \card
    \left(
    \| G_{ij} \|_F
    \right)
    $;
    $
    \left.
    \sum_{i,j}
    \log (1 + \| \gij \|_F / \veps)
    \right\},
    $
and the minimizers of~(\ref{eq.min_Gij}) are obtained by replacing the absolute value of $V_{ij}$ in~(\ref{eq.shrinkage_wl1}),~(\ref{eq.shrinkage_card}), and~(\ref{eq.rpm}) with the Frobenius norm $\| \cdot \|_F$ of the corresponding block submatrix $V_{ij}$.}

    \vspace*{-2ex}
\subsection{Anderson-Moore method for the $F$-minimization problem~{\em (\ref{eq.F_update})\/}}
    \label{sec.F_update}

We next employ the Anderson-Moore method to solve the $F$-minimization problem~(\ref{eq.F_update}). The advantage of this algorithm lies in its fast convergence (compared to the gradient method) and in its simple implementation (compared to Newton's method); e.g., see~\cite{maktoi87,rausac97,linfarjovTAC11al}. When applied to the $F$-minimization problem~(\ref{eq.F_update}), this method requires the solutions of two Lyapunov equations and one Sylvester equation in each iteration. {We next recall the first and second order derivatives of $J$; for related developments, see~\cite{rausac97}.
    \begin{proposition}
The gradient of $J$ is determined by
    \be
    \non
    \nabla J(F)
    \; = \;
    2 \, ( R F  - B_2^T P ) \, L
    \ee
where $L$ and $P$ are the controllability and observability Gramians of the closed-loop system,
    \begin{align}
      \tag*{(NC-L)}
      \label{NC1}
      ( A \,-\, B_2 F )
      \,
      L
      \,+\,
      L
      \,
      ( A \,-\, B_2 F )^T
      \,&=\,
      - \, B_1 B_1^T
      \\
      \tag*{(NC-P)}
      \label{NC2}
      ( A \,-\, B_2 F )^T \, P
      \,+\,
      P \, ( A \,- \,B_2 F )
      \,&=\,
      - \,
      ( Q  \,+\,  F^T R F ).
    \end{align}
The second-order approximation of $J$ is determined by
    \be
    \non
    J(F + \tilde{F})
    \,\approx\,
    J(F)
    \,+\,
    \langle
    \nabla
    J(F), \tilde{F}
    \rangle
    \,+\,
    (1/2)
    \,
    \langle H(F,\tilde{F}), \tilde{F} \rangle
    \ee
where $H(F,\tilde{F})$ is the linear function of $\tilde{F}$,
    \[
        H(F,\tilde{F})
        \, = \,
        2
        \left(
        (  R \tilde{F}  - B_2^T \tilde{P} ) \, L
        \,+\,
        (  R F  - B_2^T P  ) \, \tilde{L}
        \right)
    \]
and $\tilde{L}$, $\tilde{P}$ are the solutions of the following Lyapunov equations
    \[
    \ba{l}
        (A-B_2 F  )
        \,
        \tilde{L}
        \, + \,
        \tilde{L}
        \,
        ( A - B_2 F  )^T
        \, = \,
        B_2 \tilde{F}  L
        \, + \,
        (B_2 \tilde{F}  L) ^T
    \\
        (A - B_2 F  )^T \, \tilde{P}
        \, + \,
        \tilde{P} \, ( A - B_2 F  )
        \, = \,
        ( P B_2 - F^T R ) \, \tilde{F}
        \, + \,
        \tilde{F}^T  \, ( B_2^T P - R F  ).
    \ea
    \]
   \end{proposition}}

By completing the squares with respect to $F$ in the augmented Lagrangian $\cL_\rho$, we obtain the following equivalent problem to~(\ref{eq.F_update})
    \be
    \non
    \text{minimize}
    ~~
    \varphi(F)
    \, = \,
    J(F)
    \,+\,
    (\rho/2)
    \| F \,-\, U^{k} \|_F^2
    \ee
where
    $
    U^{k}
    =
    G^k
    -
    (1/\rho)
    \Lambda^k.
    $
{Setting $\nabla \varphi : = \nabla J + \rho(F - U^k)$ to zero yields the necessary conditions for optimality
    \be
        2 \left(
        R F
        \, - \,
        B_2^T P
        \right)
        L
        \, + \,
        \rho
        \left(
        F \, - \, U^{k}
        \right)
        \, = \, 0
      \tag*{(NC-F)}
      \label{NC3}
    \ee
where $L$ and $P$ are determined by~\ref{NC1} and~\ref{NC2}.}

Starting with a stabilizing feedback $F$, the Anderson-Moore method solves the two Lyapunov equations~\ref{NC1} and~\ref{NC2}, and then solves the Sylvester equation~\ref{NC3} to obtain a new feedback gain $\bar{F}$. In other words, it alternates between solving~\ref{NC1} and~\ref{NC2} for $L$ and $P$ with $F$ being fixed and solving~\ref{NC3} for $F$ with $L$ and $P$ being fixed. {It can be shown that the difference between two consecutive steps $\tilde{F} = \bar{F} - F$ forms a descent direction of $\varphi$; see~\cite{linfarjovTAC11al} for a related result.} Thus, line search methods~\cite{nocwri06} can be employed to determine step-size $s$ in $F + s \tilde{F}$ to guarantee closed-loop stability and the convergence to a stationary point of $\varphi$.

    \remark[Closed-loop stability]
    \label{rem.stab}
Since the $\htwo$ norm is well defined for causal, strictly proper, stable closed-loop systems, we set $J$ to infinity if $A - B_2 F$ is not Hurwitz. Furthermore, $J$ is a smooth function that increases to infinity as one approaches the boundary of the set of stabilizing gains~\cite{linfarjovTAC11al}. Thus, the decreasing sequence of $\{\varphi(F^i)\}$ ensures that $\{F^i\}$ are stabilizing gains.
    \vspace*{-2ex}
\subsection{Solving the structured $\htwo$ problem}
    \label{sec.polishing}

We next turn to the $\htwo$ problem subject to structural constraints on the feedback gain. Here, we fix the sparsity patterns $F \in \eS$ identified using ADMM and then solve~\ref{SH2} to obtain the optimal feedback gain that belongs to $\eS$. This procedure, commonly used in optimization~\cite[Section 6.3.2]{boyvan04}, can improve the performance of sparse feedback gains \mbox{resulting from the ADMM algorithm.}

As noted in Remark~\ref{rem.stab}, the sparse feedback gains obtained in ADMM are stabilizing. This feature  facilitates the use of descent algorithms (e.g., Newton's method) to solve~\ref{SH2}. Given an initial gain $F^0 \in \eS$, a decreasing sequence of the objective function $\{J(F^i)\}$ is generated by updating $F$ according to
    $
        F^{i+1}
         =
        F^i
         +
        s^i \, \tilde{F}^i
    $;
{here, $s^i$ is the step-size and $\tilde{F}^i \in \eS$ is the Newton direction that is determined by the minimizer of the second-order approximation of the objective function~(\ref{eq.J}). Equivalently, $\tilde{F}^i \in \eS$ is the minimizer of
    $
    \Phi(\tilde{F})
    :=
    (1/2) \,
    \langle
    H(\tilde{F}) \circ I_\eS, \tilde{F}
    \rangle
    +
    \langle
    \nabla J \circ I_\eS, \tilde{F}
    \rangle
    $
where {\em structural identity\/} $I_\eS$ of subspace $\eS$ (under entry-wise multiplication $\circ$ of two matrices) is used to characterize structural constraints
     \[
     {I_\eS}_{ij}
     \, = \,
     \left\{
     \ba{ll}
     \!\!
     1,  & \text{if}~ F_{ij} ~\text{is a free variable}
     \\[0.1cm]
     \!\!
     0,  & \text{if}~ F_{ij} = 0 ~\text{is required}
     \ea
     \right.
     \;\;\;
     \Rightarrow
     \;\;\;
     F \, \circ \, I_\eS \, = \, F
     ~~\text{for}~~
     F \, \in \, \eS.
     \]

To compute Newton direction, we use the conjugate gradient method that does not require forming or inverting the large Hessian matrix explicitly; see~\cite[Chapter 5]{nocwri06}. It is noteworthy that techniques such as the negative curvature test~\cite[Section 7.1]{nocwri06} can be employed to guarantee the descent property of the Newton direction; consequently, line search methods, such as the Armijo rule~\cite[Section 3.1]{nocwri06}, can be used to generate a decreasing sequence of $J$.}

\subsection{Convergence of ADMM}
    \label{sec.convergence}

	For convex problems the convergence of ADMM to the global minimizer follows from standard results~\cite{boyparchupeleck11}. For nonconvex problems, where convergence results are not available, extensive computational experience suggests that ADMM works well when the value of $\rho$ is sufficiently large~\cite{forglo83,glotal89}. This is attributed to the quadratic term $(\rho/2) \| F - G \|_F^2$ that tends to locally convexify the objective function for sufficiently large $\rho$; see~\cite[Chapter 14.5]{lueye08}.

	{For problem~\ref{SP} with $g$ determined by the weighted $\ell_1$ norm~(\ref{eq.wl1}), we next show that when ADMM converges, it converges to a critical point of~\ref{SP}. For a convergent point $( F^\star, G^\star, \Lambda^\star )$ of the sequence $\{ F^k,  G^k , \Lambda^k \}$, (\ref{eq.lambda_update}) simplifies to
    $
    F^\star
    -
    G^\star
    = 0.
    $
Since $F^\star$ minimizes $\cL_\rho (F,G^\star,\Lambda^\star)$ and since $G^\star$ minimizes $\cL_\rho (F^\star,G,\Lambda^\star)$, we have
   $
   \{
    0
    =
    \nabla J(F^\star)
    +
    \Lambda^\star,
    $
    $
    0
    \in
    \gamma
    \,
    \partial g(G^\star)
    -
    \Lambda^\star
    \}
    $
where $\partial g$ is the subdifferential of the convex function $g$ in~(\ref{eq.wl1}). Therefore, $(F^\star, G^\star)$ satisfies the necessary conditions for the optimality of~\ref{SP} and ADMM converges to a critical point of~\ref{SP}.}

\section{Examples}
    \label{sec.example}

We next use three examples to illustrate the utility of the approach developed in Section~\ref{sec.ADMM}. The identified sparsity structures result in {\em localized\/} controllers in all three cases. Additional information about these examples, along with {\sc Matlab} source codes, can be found at
    \vspace*{1ex}
    \begin{center}
    \href{http://www.ece.umn.edu/~mihailo/software/lqrsp/index.html}{\sf www.ece.umn.edu/$\sim$mihailo/software/lqrsp/}
    \end{center}

    \vspace*{-2ex}
\subsection{Mass-spring system}
    \label{sec.mass-spring}

For a mass-spring system with $N$ masses on a line, let $p_i$ be the displacement of the $i$th mass from its reference position and let the state variables be $x_1 :=  [\, p_1 \, \cdots \, p_N \,]^T$ and $x_2 := \dot{x}_1$. For unit masses and spring constants, the state-space representation is given by~(\ref{eq.ss}) with
    \be
    \non
     A
     \; = \;
     \left[
     \ba{cc}
     O & I\\
     T & O
      \ea
     \right],
    ~~~
    B_1
    \; = \;
    B_2
    \; = \;
    \left[
    \ba{cc}
     O \\
     I
     \ea
     \right] ,
    \ee
where $T$ is an $N \times N$ tridiagonal Toeplitz matrix with $-2$ on its main diagonal and $1$ on its first sub- and super-diagonal, and $I$ and $O$ are $N \times N$ identity and zero matrices. The state performance weight $Q$ is the identity matrix and the control performance weight is $R = 10I$.

We use the weighted $\ell_1$ norm as the sparsity-promoting penalty function, where we follow~\cite{canwakboy08} and set the weights $W_{ij}$ to be inversely proportional to the magnitude of the solution $F^\star$ of~\ref{SP} at the previous value of $\gamma$,
    $
    W_{ij}
    =
    1/(
    |F_{ij}^\star|
    +
    \veps).
    $
This places larger relative weight on smaller feedback gains and they are more likely to be dropped in the sparsity-promoting algorithm. Here, $\veps = 10^{-3}$ is introduced to have well-defined weights when $F^\star_{ij} = 0$.

\begin{figure}
      \begin{tabular}{cccc}
      \!\!\!\!\!\!\!\!
      \begin{tabular}{c}
      \subfloat[$\gamma = 10^{-4}$]
      {\label{fig.ms_Fp}
      \includegraphics[width=0.24\textwidth]{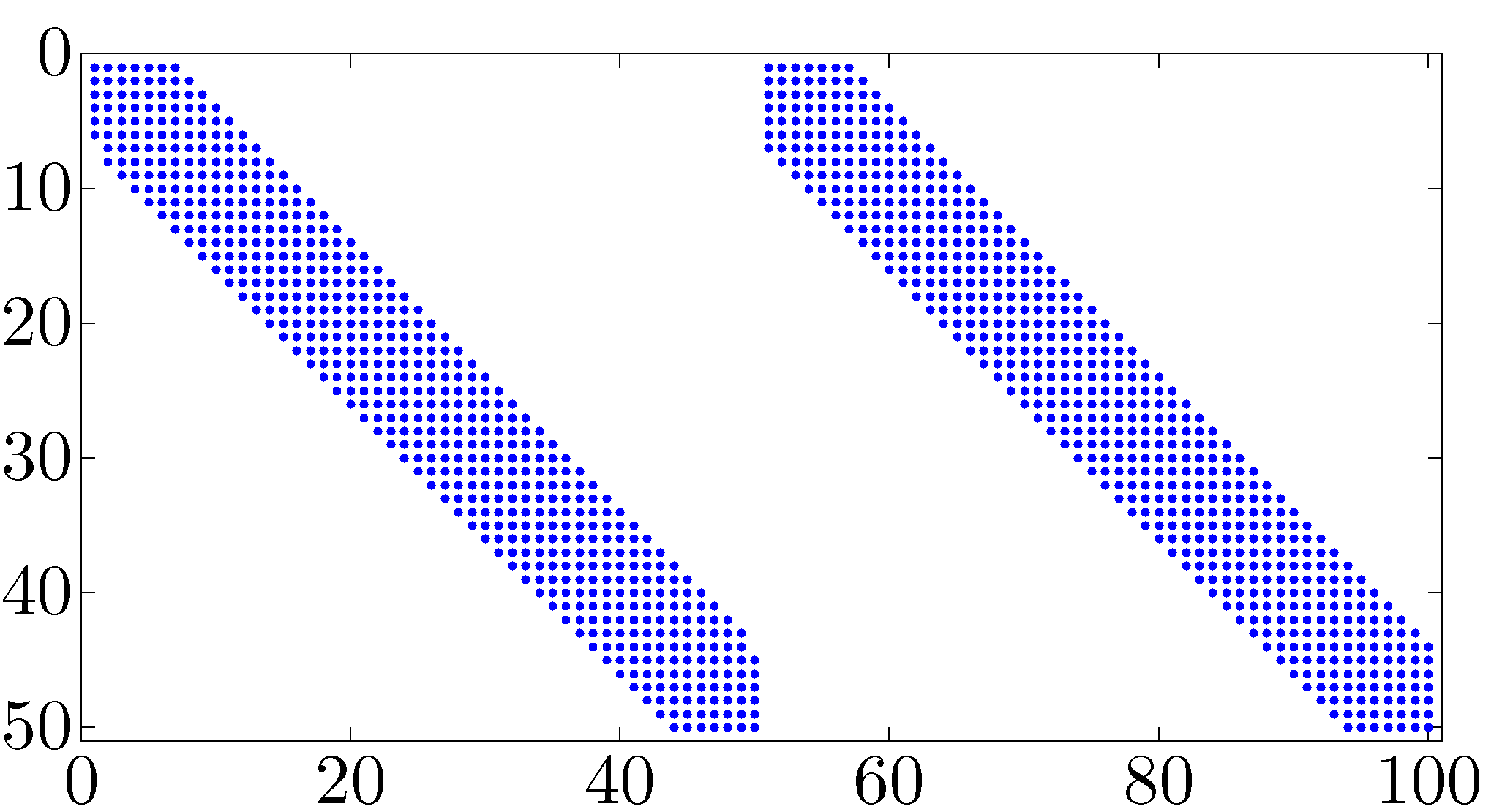}}
      \end{tabular}
      \!\!\!\!\!\!\!\!&\!\!\!\!\!\!\!\!
      \begin{tabular}{c}
      \subfloat[$\gamma = 0.0105$]
      {\label{fig.ms_Fv}
      \includegraphics[width=0.24\textwidth]{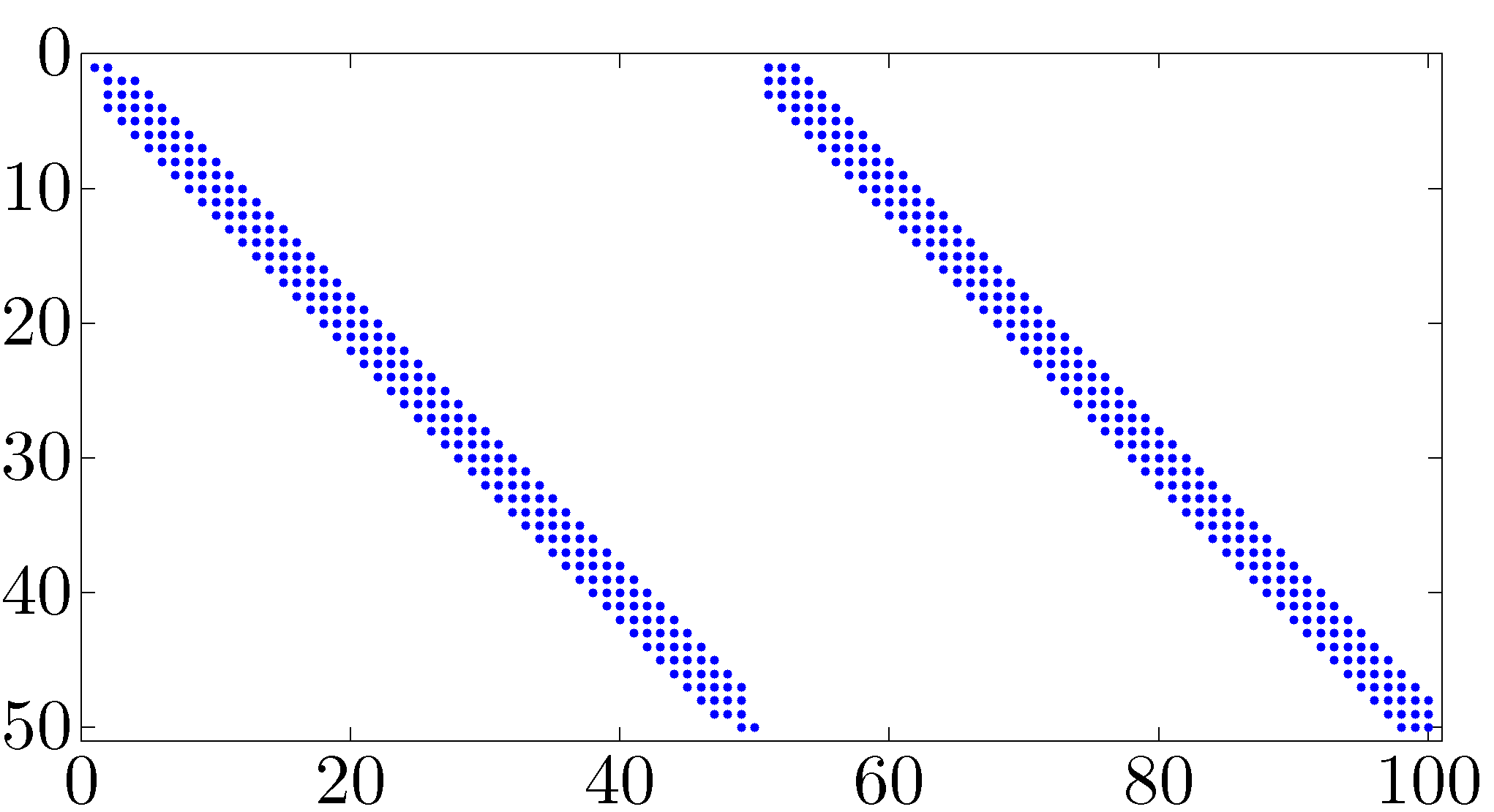}}
      \end{tabular}
      \!\!\!\!\!\!\!\!&\!\!\!\!\!\!\!\!
      \begin{tabular}{c}
      \subfloat[]
      {\label{fig.ms_diag_Fp}
      \includegraphics[width=0.24\textwidth]{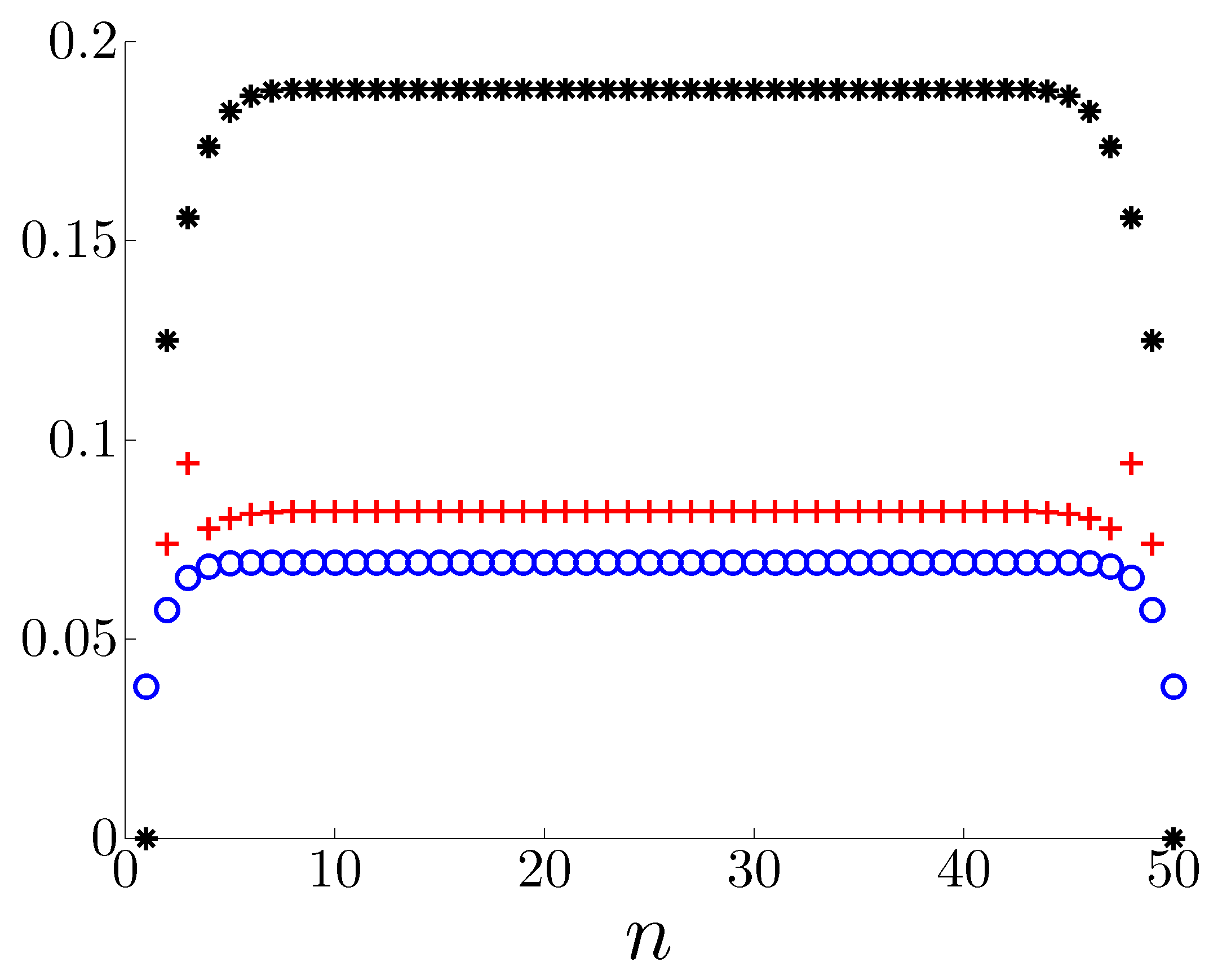}}
      \end{tabular}
      \!\!\!\!\!\!\!\!&\!\!\!\!\!\!\!\!
      \begin{tabular}{c}
      \subfloat[]
      {\label{fig.ms_diag_Fv}
      \includegraphics[width=0.24\textwidth]{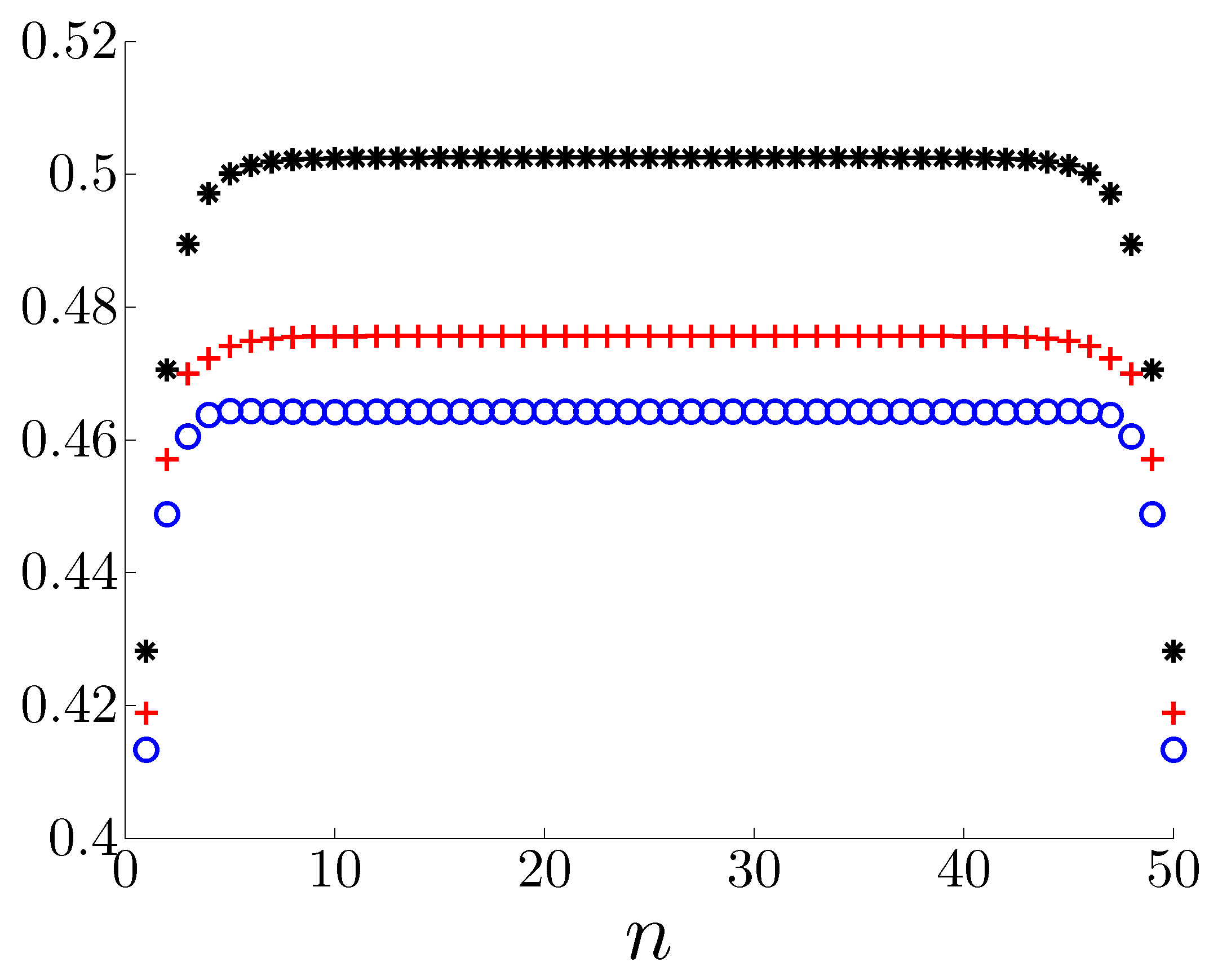}}
      \end{tabular}
      \end{tabular}
      \caption{Sparsity patterns of $F^\star = [F^\star_p~~F^\star_v] \in \mathbb{R}^{50 \times 100}$ for the mass-spring system obtained using weighted $\ell_1$ norm with (a) $\gamma = 10^{-4}$ and (b) $\gamma = 0.0105$. As $\gamma$ increases, the number of nonzero sub- and super-diagonals of $F_p^\star$ and $F_v^\star$ decreases. The diagonals of (c) $F_p^\star$ and (d) $F_v^\star$ for different values of $\gamma$: $10^{-4}$ ($\circ$), $0.0281$ ({\small $+$}), and $0.1$ ({\small $*$}). The diagonals of the centralized position and velocity gains are almost identical to ($\circ$).}
      \label{fig.ms_F}
    \end{figure}

The optimal feedback gain at $\gamma = 0$ is computed from the solution of the algebraic Riccati equation. As $\gamma$ increases, the number of nonzero sub- and super-diagonals of both position $F^\star_p$ and velocity $F^\star_v$ gains decreases; see {Figs.~\ref{fig.ms_Fp} and~\ref{fig.ms_Fv}.} Eventually, both $F^\star_p$ and $F^\star_v$ become diagonal matrices. It is noteworthy that diagonals of both position and velocity feedback gains are nearly constant except for masses that are close to the boundary; see~Figs.~\ref{fig.ms_diag_Fp} and~\ref{fig.ms_diag_Fv}.

\begin{figure}
      \centering
      \begin{tabular}{ccc}
      \!\!\!\!\!\!\!\!\!\!\!
      \begin{tabular}{c}
      \subfloat[${\bf card} \, (F^\star) / {\bf card} \, (F_{c})$]
      {
      \includegraphics[width=0.24\textwidth]{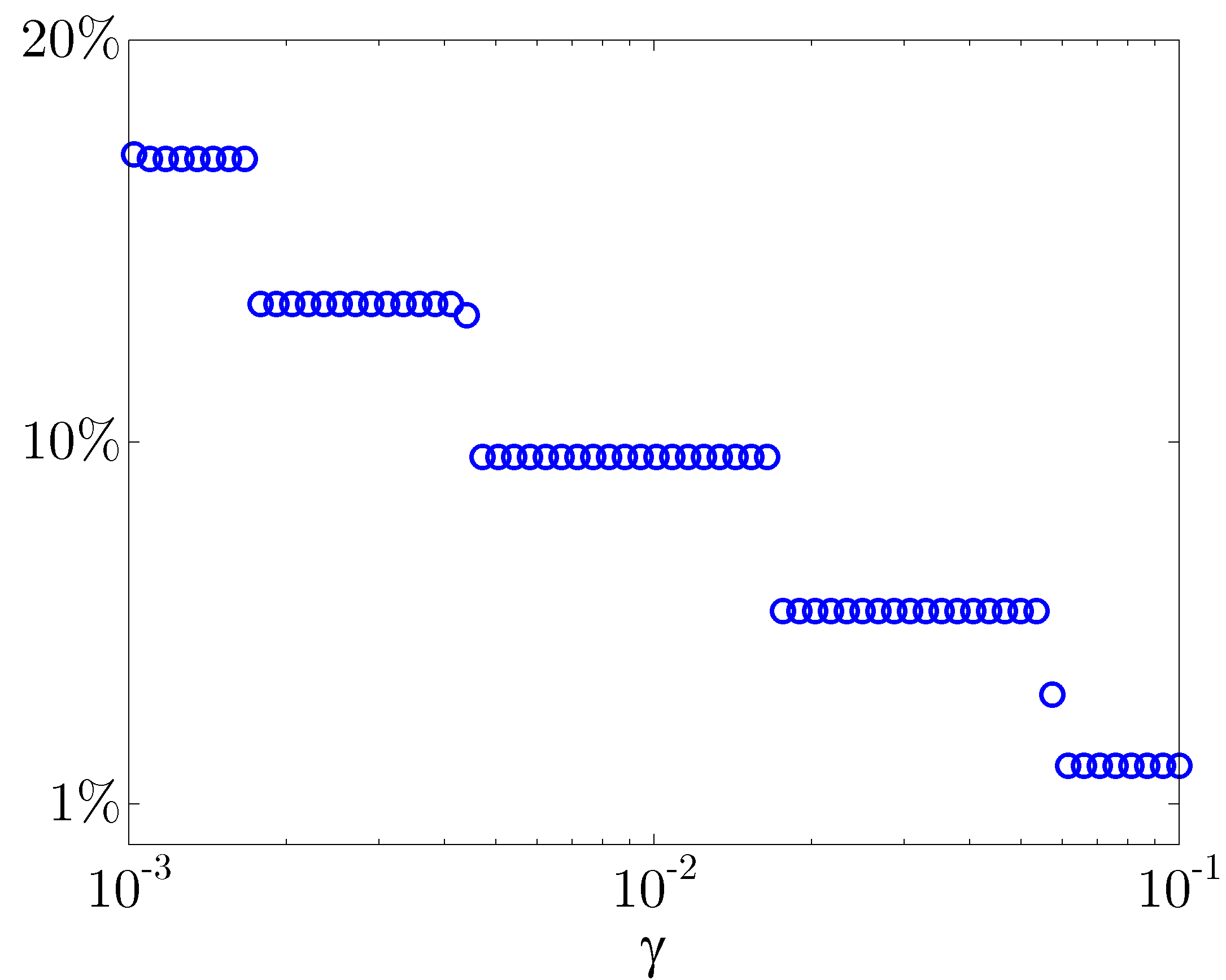}}
      \end{tabular}
      \!\!\!\!\!\! & \!\!\!\!\!\!\!\!
      \begin{tabular}{c}
      \subfloat[$(J(F^\star) \,-\, J(F_c)) / J(F_c)$]
      {
      \includegraphics[width=0.24\textwidth]{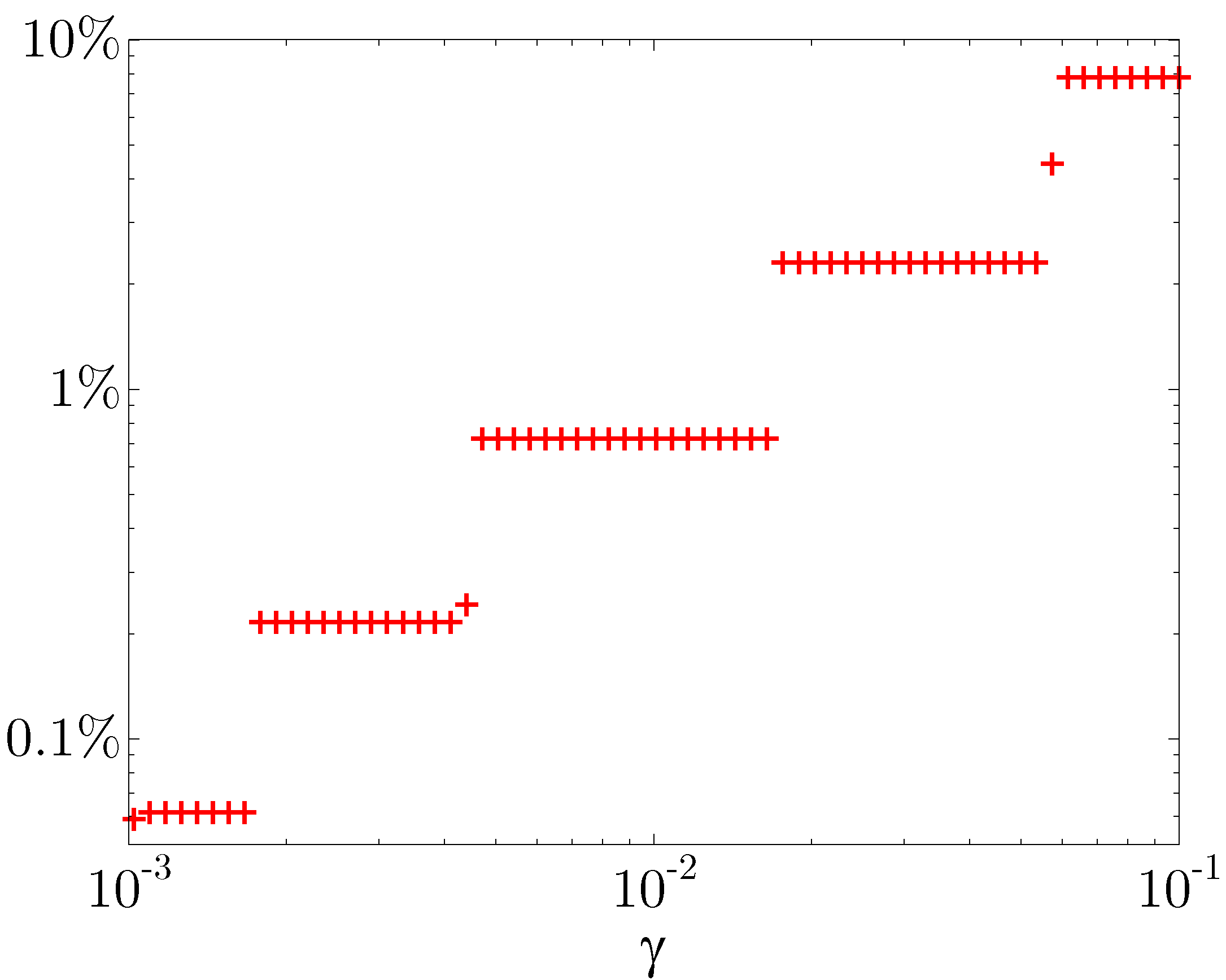}}
      \end{tabular}
      \!\!\!\!\!\!\!\! & \!\!\!\!\!\!
      \subfloat[]{
        \begin{tabular}{clll}
        $\gamma$
        &
        $0.01$
        &
        $0.04$
        &
        $0.10$
        \\[0.1cm]
        ${\bf card} \, (F^\star) / {\bf card} \, (F_{c})$
        &
        $9.4 \%$
        &
        $5.8 \%$
        &
        $2.0 \%$
        \\[0.1cm]
        $(J(F^\star) \,-\, J(F_c)) / J(F_c)$
        &
        $0.8 \%$
        &
        $2.3 \%$
        &
        $7.8 \%$
        \end{tabular}}
      \end{tabular}
      \caption{(a) The sparsity level and (b) the performance degradation of $F^\star$ compared to the centralized gain $F_c$ for mass-spring system. (c) Sparsity vs.\ performance: using $2\%$ of nonzero elements, $\htwo$ performance of $F^\star$ is only $7.8\%$ worse than performance of $F_c$.}
      \label{fig.ms_sparsity}
    \end{figure}

After sparsity structures of controllers are identified by solving~\ref{SP}, we fix sparsity patterns and solve structured $\htwo$ problem~\ref{SH2} to obtain the optimal {\em structured\/} controllers. Comparing the sparsity level and the performance of these controllers to those of the centralized controller $F_c$, we see that using only a {\em fraction\/} of nonzero elements, the sparse feedback gain $F^\star$ achieves $\htwo$ performance comparable to the performance of $F_c$; see Fig.~\ref{fig.ms_sparsity}. In particular, using about $2\%$ of nonzero elements, $\htwo$ performance of $F^\star$ is only about $8\%$ worse than that of $F_c$.

    \vspace*{-2ex}
\subsection{Network with $100$ unstable nodes}

Let $N = 100$ nodes be randomly distributed with a uniform distribution in a square region of $10 \times 10$ units. Each node is an unstable second order system coupled with other nodes through the exponentially decaying function of the Euclidean distance $\alpha(i,j)$ between them~\cite{motjad08}
    \be
    \non
    \tbo{\dot{x}_{1i}}{\dot{x}_{2i}}
    \, = \,
    \left[
    \begin{array}{cc}
    1 & 1 \\
    1 & 2
    \end{array}
    \right]
    \tbo{x_{1i}}{x_{2i}}
    \; + \;
    \sum_{
    j \,\neq\, i}
    \mre^{- \alpha(i,j)}
    \tbo{x_{1j}}{x_{2j}}
    \; + \;
    \left[
    \begin{array}{c}
    0 \\
    1
    \end{array}
    \right]
    \left(
    d_i
    \,+\,
    u_i
    \right)
    \ee
with $i = 1, \ldots, N$. The performance weights $Q$ and $R$ are set to identity matrices.

\begin{figure}
      \centering
      \subfloat[$\gamma=12.6$]
      {\label{fig.comm_graph1}
      \includegraphics[width=0.22\textwidth]{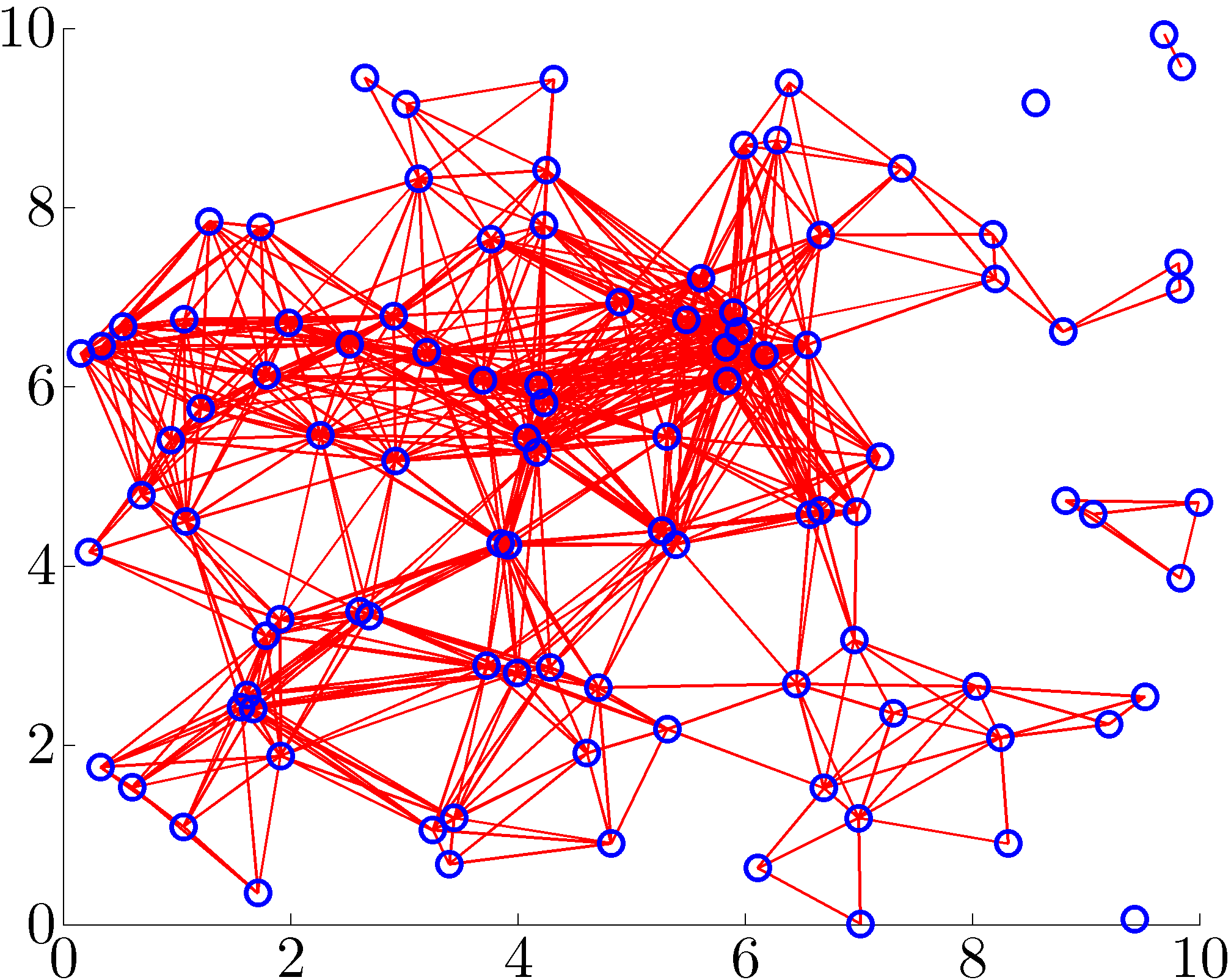}}
      ~
      \subfloat[$\gamma=26.8$]
      {\label{fig.comm_graph2}
      \includegraphics[width=0.22\textwidth]{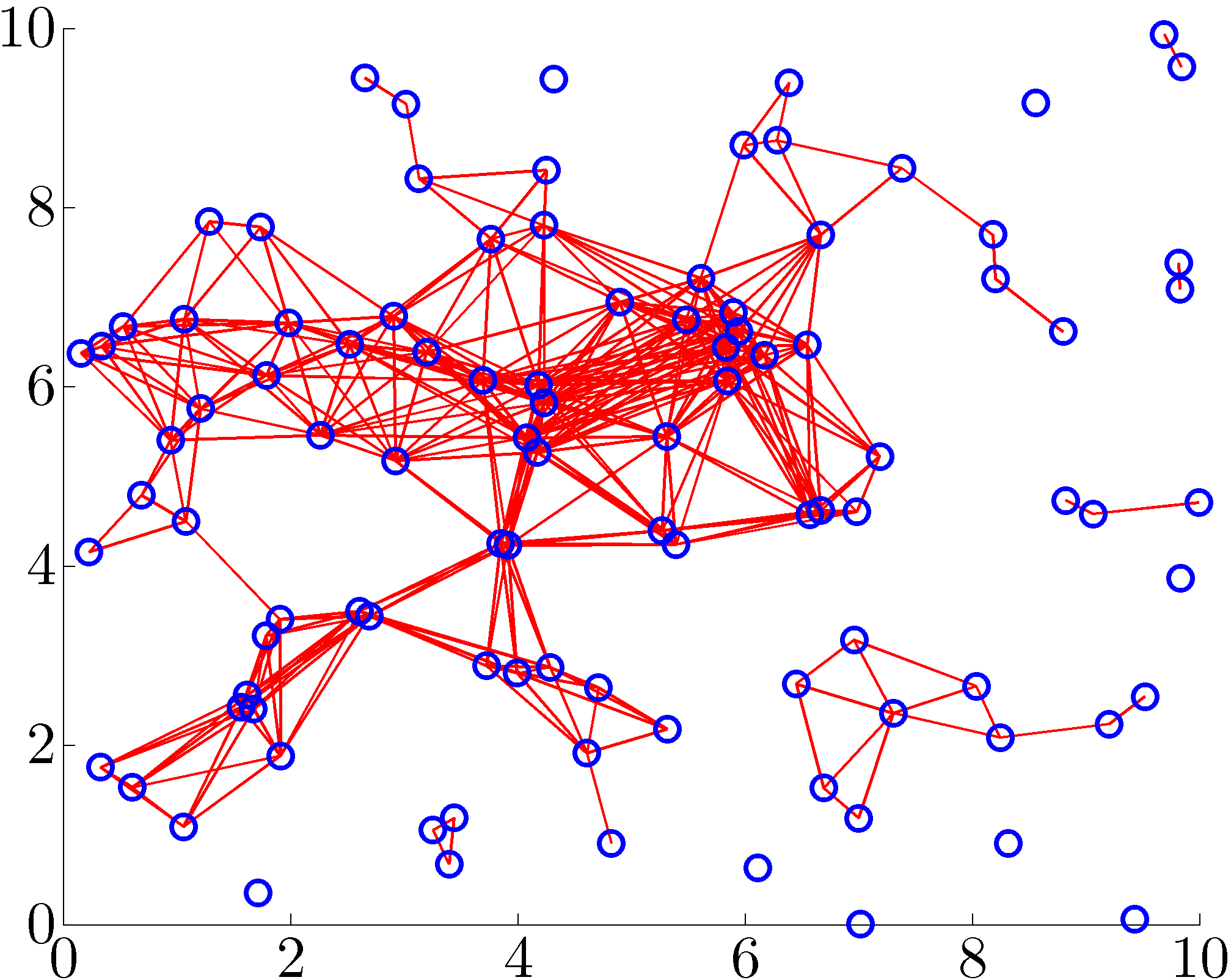}}
      ~
      \subfloat[$\gamma=68.7$]
      {\label{fig.comm_graph3}
      \includegraphics[width=0.22\textwidth]{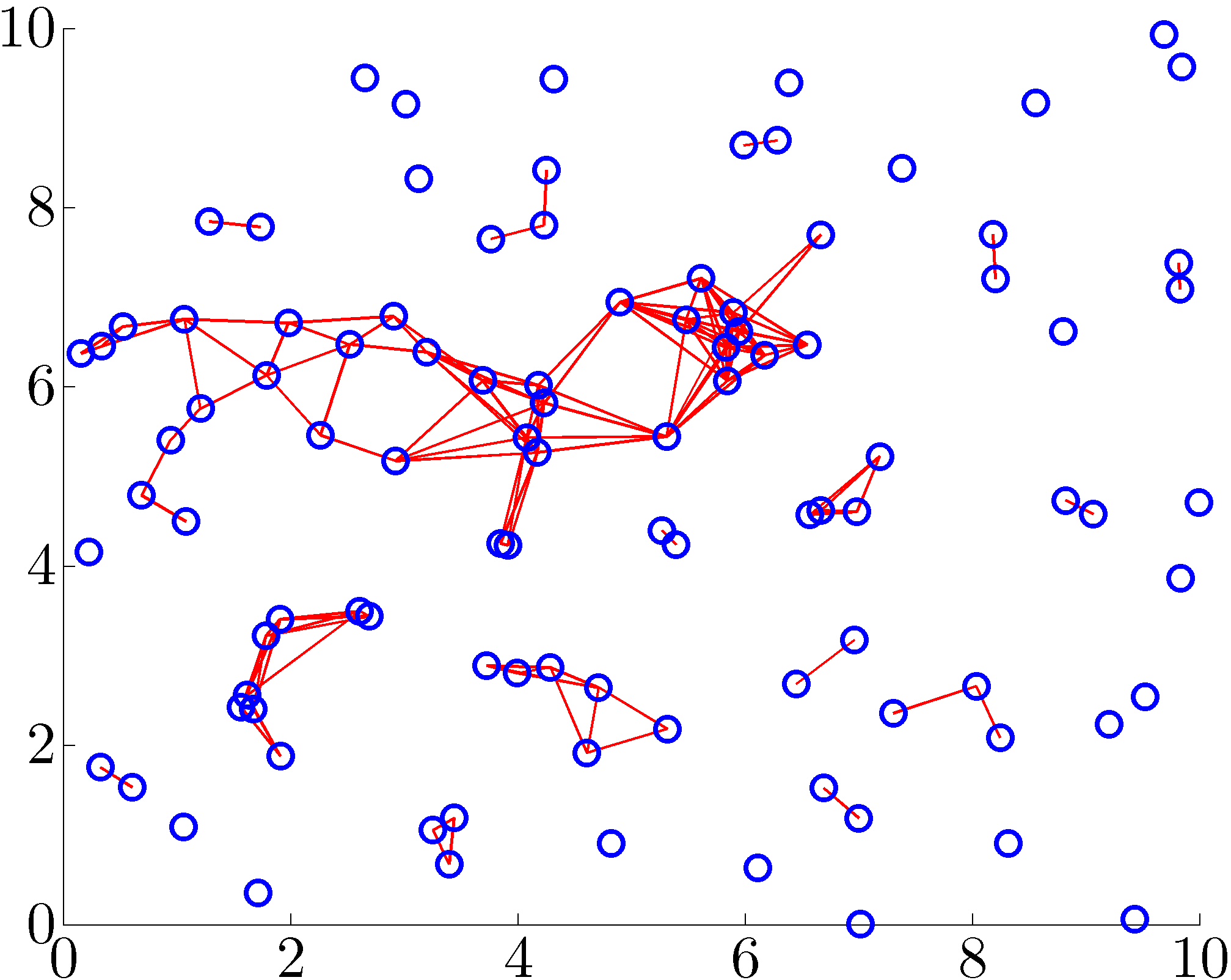}}
      ~
      \subfloat[]
      {\label{fig.H2_vs_sp}
      \includegraphics[width=0.24\textwidth]{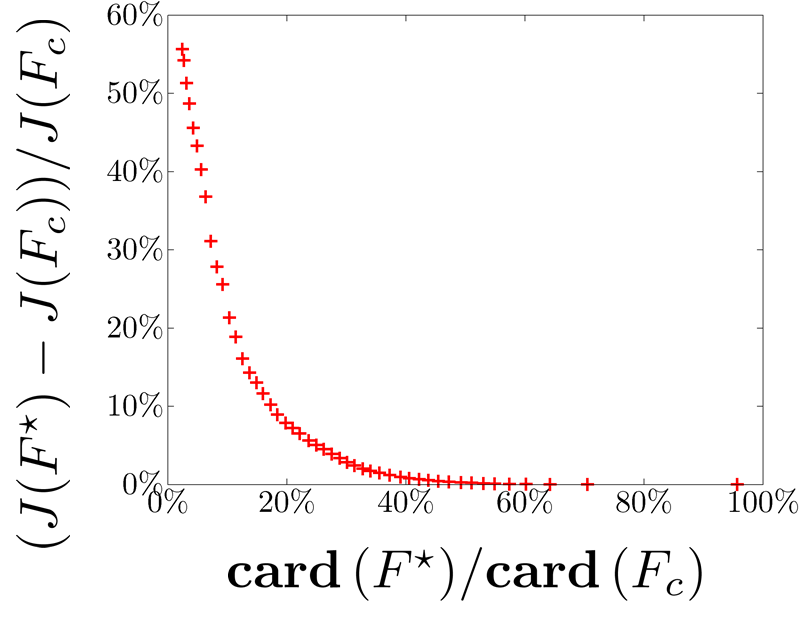}}
      \caption{(a)-(c) The localized communication graphs of distributed controllers obtained by solving~\ref{SP} for different values of $\gamma$  for the network with $100$ nodes. Note that the communication graph does not have to be connected since the nodes are dynamically coupled to each other and allowed to measure their own states. (d) The optimal trade-off curve between the $\htwo$ performance degradation and the sparsity level of $F^\star$ compared to the centralized gain $F_c$.}
      \label{fig.comm_graph}
    \end{figure}

We use the weighted $\ell_1$ norm as the penalty function with the weights given in Section~\ref{sec.mass-spring}. As $\gamma$ increases, the underlying communication graphs gradually become {\em localized\/}; see Figs.~\ref{fig.comm_graph1},~\ref{fig.comm_graph2}, and~\ref{fig.comm_graph3}. With about $8\%$ of nonzero elements of $F_c$, $\htwo$ performance of $F^\star$ is about $28\%$ worse than performance of the centralized gain $F_c$. Figure~\ref{fig.H2_vs_sp} shows the optimal trade-off curve between the $\htwo$ performance and the feedback gain sparsity.

We note that the truncation of the centralized controller could result in a non-stabilizing feedback matrix~\cite{motjad08}. In contrast, our approach gradually modifies the feedback gain and increases the number of zero elements, which plays an important role in preserving the closed-loop stability.

    \vspace*{-2ex}
\subsection{Block sparsity: A bio-chemical reaction example}

Consider a network of $N = 5$ systems coupled through the following dynamics
    \[
    \dot{x}_i
    \, = \,
    [A]_{ii}
    \,
    x_i
    \, - \,
    \frac{1}{2}
    \sum_{
    j \,=\, 1
    }^N
    (i -  j)
    \,
    (x_i
    \, - \,
    x_j)
    \, + \,
    [B_1]_{ii}
    \,
    d_i
    \,+\,
    [B_2]_{ii}
    \,
    u_i
    \]
where $[\, \cdot \,]_{ij}$ denotes the $ij$th block of a matrix and
    \[
    [A]_{ii}
    \, = \,
    \left[
    \begin{array}{rrr}
    -1  & 0  & -3 \\
    3   & -1 & 0    \\
    0   & 3  & -1
    \end{array}
    \right],
    ~~
    [B_1]_{ii}
    \, = \,
    \left[
    \begin{array}{ccc}
    3   & 0  & 0  \\
    0   & 3  & 0  \\
    0   & 0  & 3
    \end{array}
    \right],
    ~~
    [B_2]_{ii}
    \, = \,
    \left[
    \begin{array}{c}
    3 \\
    0 \\
    0
    \end{array}
    \right].
    \]
The performance weights $Q$ and $R$ are set to identity matrices. Systems of this form arise in bio-chemical reactions with a cyclic negative feedback~\cite{jovarcsonTAC08}.

\begin{figure}
        \centering
        \begin{tabular}{cc}
        \subfloat[]
        {\label{fig.5nodes_F_block}
        \includegraphics[width=0.45\textwidth]{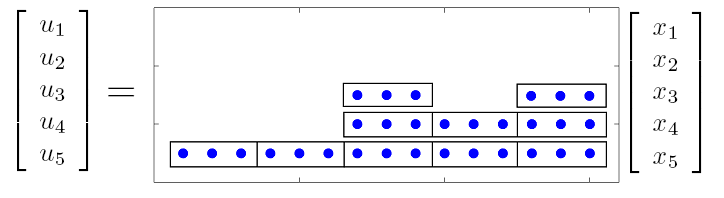}}
        &
        \subfloat[]
        {\label{fig.5nodes_F}
        \includegraphics[width=0.45\textwidth]{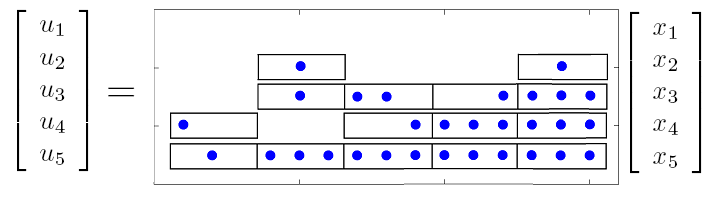}}
        \\[-0.4cm]
        \subfloat[]
        {\label{fig.5nodes_comm_F_block}
        \includegraphics[width=0.2\textwidth]{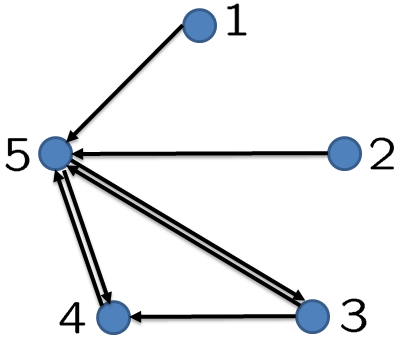}}
        &
        \subfloat[]
        {\label{fig.5nodes_comm_F}
        \includegraphics[width=0.2\textwidth]{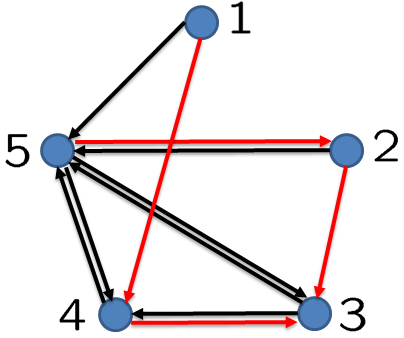}}
        \end{tabular}
        \caption{The sparse feedback gains obtained by solving~\ref{SP} using (a) the weighted sum of Frobenius norms with $\gamma = 3.6$ and (b) the weighted $\ell_1$ norm~(\ref{eq.wl1}) with $\gamma = 1.3$. Here, $F \in \mathbb{R}^{5 \times 15}$ is partitioned into $25$ blocks $F_{ij} \in \mathbb{R}^{1 \times 3}$. Both feedback gains have the same number of nonzero elements (indicated by dots) and close $\htwo$ performance (less than $1\%$ difference), but different number of nonzero blocks (indicated by boxes). Communication graphs of (c) the block sparse feedback gain in (a), and (d) the sparse feedback gain in (b) (red color highlights the additional links). An arrow pointing from node $i$ to node $j$ indicates that $i$ uses measurements from $j$.}
        \label{eq.5nodes_gains}
    \end{figure}

We use the weighted sum of Frobenius norms as the sparsity-promoting penalty function and we set the weights $W_{ij}$ to be inversely proportional to the Frobenius norm of the solution $F^\star_{ij}$ to~\ref{SP} at the previous value of $\gamma$, i.e.,
    $
    W_{ij}
    =
    1/
    (
    \| F_{ij}^\star \|_F
    +
    \veps)
    $
with
    $
    \veps
    =
    10^{-3}.
    $
As $\gamma$ increases, the number of {\em nonzero blocks\/} in $F$ decreases. Figure~\ref{eq.5nodes_gains} shows sparsity patterns of feedback gains and the corresponding communication graphs resulting from solving~\ref{SP} with {sparse and block sparse penalty functions.} Setting $\gamma$ to values that yield the same number of nonzero elements in these feedback gains results in the block sparse feedback gain with a smaller number of nonzero blocks. In particular, the first two rows of the block sparse feedback gain in Fig.~\ref{fig.5nodes_F_block} are identically equal to zero (indicated by blank space). This means that the subsystems $1$ and $2$ do not need to be actuated. Thus, the communication graph determined by the block sparse feedback gain has fewer links; cf.\ Figs.~\ref{fig.5nodes_comm_F_block} and~\ref{fig.5nodes_comm_F}.

    \vspace*{-2ex}
\section{Concluding Remarks}
    \label{sec.conclusion}

We design sparse and block sparse state feedback gains that optimize the $\htwo$ performance of distributed systems. The design procedure consists of a structure identification step and a ``polishing'' step. In the identification step, we employ the ADMM algorithm to solve the sparsity-promoting optimal control problem, whose solution gradually moves from the centralized gain to the sparse gain of interest as our emphasis on the sparsity-promoting penalty term is increased. In the polishing step, we use Newton's method in conjunction with a conjugate gradient scheme to solve the minimum variance problem subject to the identified sparsity constraints.

Although we focus on the $\htwo$ performance, the developed framework can be extended to design problems with other performance indices. We emphasize that the analytical solutions to the $G$-minimization problem are independent of the assigned performance index. Consequently, the $G$-minimization step in ADMM for~\ref{SP} with alternative performance indices can be done exactly as in Section~\ref{sec.G_update}. Thus, ADMM provides a flexible framework for sparsity-promoting optimal control problems of the form~\ref{SP}.

We have recently employed ADMM for selection of an {\em a priori\/} specified number of leaders in order to minimize the variance of stochastically forced dynamic networks~\cite{linfarjovCDC11}, {for creation of new social links to maximize public awareness in social networks~\cite{farzhalinjovCDC12}, and for identification of sparse representations of consensus networks~\cite{dhilinfarjovIFAC12}.} We also aim to extend the developed framework to the observer-based sparse optimal feedback design. Our results on the identification of classes of convex optimal control problems will be reported elsewhere.

\vspace*{-2ex}
\section*{Acknowledgements}

We would like to thank Stephen P.\ Boyd for inspiring discussions on ADMM, group lasso, and cardinality minimization, and Roland Glowinski for his insightful comments on convergence of ADMM.

\end{document}